\documentclass[10pt,onecolumn]{IEEEtran}
\IEEEoverridecommandlockouts

\usepackage{amsbsy}
\usepackage{amsmath}
\usepackage{amssymb}
\usepackage{latexsym}
\usepackage{ifthen}
\usepackage{subfigure}
\usepackage{turnstile}
\usepackage{mathtools}
\usepackage{booktabs}
\usepackage{balance}
\usepackage{paralist}
\usepackage{cite}
\usepackage{url}
\usepackage{color}

\newtheorem{theorem}{Theorem}[section]

\newcommand{\spr}{\boldsymbol{\operatorname{spr}}}

\newcommand{\meas}{\boldsymbol{\operatorname{meas}}}
\newcommand{\cl}{\boldsymbol{\operatorname{cl}}}

\def\mb{\mathbf}

\def\z{\mb{z}}

\def\u{\mb{u}}

\def\v{\mb{v}}

\def\oA{\overline{A}}

\def\v{\mb{v}}

\def\u{\mb{u}}

\def\mb{\mathbf}

\newtheorem{definition}{\bf Definition}

\def \BNN{\{0,1\}^{N\times N}}

\newcommand{\cmat}[2]{\mathbb{C}^{#1\times #2}}

\title{Finite-Time In-Network Computation of Linear Transforms}

\author{Soummya Kar,\IEEEauthorrefmark{2} Markus P\"uschel,\IEEEauthorrefmark{3} and Jos\'e M. F. Moura\IEEEauthorrefmark{2}\thanks{This work was supported in part by NSF under grant CCF-1513936.}}

\begin{document}

\maketitle

\renewcommand{\thefootnote}{\fnsymbol{footnote}}

\footnotetext[2]{Department of Electrical and Computer Engineering, Carnegie Mellon University, Pittsburgh, PA, USA (soummyak@andrew.cmu.edu, moura@andrew.cmu.edu).}
\footnotetext[3]{Department of Computer Science, ETH Zurich, Z\"urich, Switzerland (pueschel@inf.ethz.ch).}

\renewcommand{\thefootnote}{\arabic{footnote}}

\thispagestyle{empty}


%


\maketitle

\begin{abstract} This paper focuses on finite-time in-network computation of linear transforms of distributed graph data. Finite-time transform computation problems are of interest in graph-based computing and signal processing applications in which the objective is to compute, by means of distributed iterative methods, various (linear) transforms of the data distributed at the agents or nodes of the graph. While finite-time computation of consensus-type or more generally rank-one transforms have been studied, systematic approaches toward scalable computing of general linear transforms, specifically in the case of heterogeneous agent objectives in which each agent is interested in obtaining a different linear combination of the network data, are relatively less explored. In this paper, by employing ideas from algebraic geometry, we develop a systematic characterization of linear transforms that are amenable to distributed in-network computation in finite-time using linear iterations. Further, we consider the general case of directed inter-agent communication graphs. Specifically, it is shown that \emph{almost all} linear transformations of data distributed on the nodes of a digraph containing a Hamiltonian cycle may be computed using at most $N$ linear distributed iterations. Finally, by studying an associated matrix factorization based reformulation of the transform computation problem, we obtain, as a by-product, certain results and characterizations on sparsity-constrained matrix factorization that are of independent interest.
\end{abstract}


\thispagestyle{plain}
\markboth{}{}

\section{Introduction}
\label{sec:intro} 
Suppose we have $N$ agents with (possibly complex-valued) data $z_{n}\in\mathbb{C}$, $i=1,\cdots,N$, distributed across the agents. Given a matrix $T\in\mathbb{C}^{N\times N}$, we are interested in computing the linear transformation
\begin{align}
\label{lin_trans}
\mathbf{y}=T\mathbf{z},
\end{align}
$\mathbf{z}=[z_{1}~z_{2}\cdots z_{N}]^{\top}$, in a \emph{distributed fashion} such that agent $n$ obtains the $n$-th component $y_{n}=\mathbf{t}_{n}^{\top}\mathbf{z}$ of the transform vector, with $\mathbf{t}_{n}^{\top}$ denoting the $n$-th row of $T$. This basic problem of distributed transform computation forms a key building block in numerous applications in distributed information processing and control ranging from distributed inference, learning, to optimization in networked multi-agent systems and control in networked dynamical systems such as cyber-physical systems. Note that, to obtain $y_{n}=\mathbf{t}_{n}^{\top}\mathbf{z}$ at agent $n$, it needs access (possibly indirectly) to data $z_{l}$ at all agents $l$ such that the corresponding entry in the transform vector $\mathbf{t}_{n}^{\top}$ is non-zero. In what follows, we are interested in distributed computation, in that, an agent may not directly communicate or exchange information with all other agents but inter-agent communication is restricted to a preassigned digraph $\mathcal{G}$ defined on the set of $N$ agents. Denote by $\bar{A}\in\{0,1\}^{N\times N}$ the (directed) unweighted adjacency matrix induced by $\mathcal{G}$, i.e., agent $l$ can directly send information to agent $n$ iff $\bar{A}_{nl}=1$. (We assume $\bar{A}_{nn}=1$ for all $n$, i.e., agents can access their own information.) Adopting the widely studied framework of consensus or gossip type distributed algorithms~\cite{dimakis2010gossip}, we propose to develop iterative approaches of the form described below to compute the linear transform $T\mathbf{z}$ of the agent data in a distributed fashion: each agent $n$ iteratively updates an information state $z_{n}(k)\in\mathbb{C}$ as a linear combination of its current information state and those of its \emph{neighbors} $l$, i.e.,
\begin{align}
\label{lin_dist_update}
z_{n}(k+1)=\sum_{l\in\Omega_{n}}a_{nl}(k)z_{l}(k),~~~t\geq 0,
\end{align}
where $k$ denotes the iteration or discrete time counter, $\Omega_{n}$ denotes the communication neighborhood of agent $n$, i.e., $l\in\Omega_{n}$ iff $\bar{A}_{nl}=1$ or agent $l$ can send its information state $z_{l}(k)$ directly to agent $n$, and the $a_{nl}(k)$'s are time-varying combination weights to be designed. The agent iterates are initialized as $z_{n}(0)=z_{n}$ for all $n$. First, note that the iterative update in~\eqref{lin_dist_update} is fully distributed in that each agent uses its local information and information in its (direct) communication neighborhood to implement its state update. Depending on the transform $T$ of interest, the goal is to design the weights $a_{nl}(k)$'s in~\eqref{lin_dist_update} such that agents get to successfully compute the transformation~\eqref{lin_trans} in a finite number $d$ of iterations, i.e., such that $z_{n}(d)=y_{n}=\mathbf{t}_{n}^{\top}\mathbf{z}$ for all $n$. In this paper we address the following question: under what conditions on the graph adjacency structure $\bar{A}$ it is possible to design the weights $a_{nl}(k)$'s in~\eqref{lin_dist_update} to achieve finite-time distributed computation of \emph{almost all}\footnote{By \emph{almost all}, we mean all transforms or matrices in $\mathbb{C}^{N\times N}$ except a set a set of matrices whose closure is nowhere dense in $\mathbb{C}^{N\times N}$. In particular, if we are able to implement almost all transforms, we can implement any transform with arbitrary accuracy.} linear transforms, and, if so, how many iterations are required to achieve such computation? The main challenge here is that the transform $T$ of interest and the inter-agent communication adjacency structure $\bar{A}$ could be very different, for instance, $\bar{A}$ might correspond to a sparsely connected digraph whereas $T$ could be a full or very dense matrix.

While the literature on asymptotic distributed linear transform computation, primarily based on consensus or gossip type procedures~\cite{dimakis2010gossip,kar2013consensus,teke2020random}, is extensive, the development in finite-time computation is relatively sparse. Most of the development so far in finite-time distributed transform computation has focused on the specific case of rank-one transforms where $T=\u\v^{\top}$ with $\u,\v\in\mathbb{C}^{N}$, i.e., in other words, the scenario in which the agents are interested in obtaining the same linear combination $\v^{\top}\z$ of the data $\z$ (up to scalings). A further special case is the celebrated average consensus problem in which $T=(1/N)\mathbf{1}\mathbf{1}^{\top}$, i.e., all agents have the homogeneous objective of computing the average of the network data. For connected undirected graphs, the rank-one case is known to be solvable within $N$ steps, where $N$ denotes the number of agents, see, for example,~\cite{sandryhaila2014finite,safavi2014revisiting,segarra2017optimal} for representative work. 

We provide a brief outline of the rest of the paper. In Section~\ref{sec:prob_form} we formally present the finite-time distributed transform computation problem and obtain a reformulation based on a graph-structured matrix factorization problem. Section~\ref{sec:fact} provides a general graph-structured matrix factorization result that is of independent interest and used to obtain solutions to the finite-time distributed transform computation problem in Section~\ref{sec:main _res}. Finally, Section~\ref{sec:disc} concludes the paper and presents avenues for future research.

\section{Problem Formulation and Preliminaries}
\label{sec:prob_form}

In this section, we set up some notation and introduce the distributed transform computation problem formally. We provide a brief review of the literature relevant to finite-time in-network computation of linear transforms. Additionally, we show that the distributed transform computation problem may be reformulated as a graph-structured matrix factorization problem that is of independent interest and forms the technical basis of our solution approach.

\noindent\textbf{Notation}. For a digraph $G=(V,E)$ on $N$ agents, $V=[1,\cdots,N]$ denotes the set of (labeled) agents and $E$ the set of edges or direct communication links among agents. The (unweighted) adjacency matrix $\oA\in\{0,1\}^{N\times N}$ of $G$ is given by, $\oA_{(n,l)}=1$ if $(n,l)\in E$, or $0$ otherwise. Unless stated otherwise, the digraphs to be studied in this paper possess self-loops at all vertices, i.e., $(n,n)\in E$ for all $n\in V$ or equivalently the associated adjacency matrix $\oA$ contains a strictly non-zero diagonal.

\noindent By $I_{r}$ we denote the $r\times r$ identity matrix.

\noindent For a matrix $P\in\cmat{N}{N}$, we will denote by $\overline{P}\in\{0,1\}^{N\times N}$ the structural matrix of $P$, i.e., $\overline{P}_{(n,l)} = 1$ if $P_{(n,l)}\neq 0$, and $\overline{P}_{(n,l)}=0$ otherwise. Thus, to any $P\in\cmat{N}{N}$, we can associate a digraph on $V$ specified by the adjacency matrix $\overline{P}\in\{0,1\}^{N\times N}$. This correspondence between digraphs on $V$ and matrices in $\cmat{N}{N}$ will be used throughout. Finally, by $\preceq$ we will denote the partial order on $\{0,1\}^{N\times N}$, such that $\overline{P}\preceq\overline{Q}$ implies $\overline{P}_{(n,l)}=0$ if $\overline{Q}_{(n,l)}=0$ for all $(n,l)\in V^{2}$.

\noindent By $\meas(\cdot)$, we denote the Lebesgue measure on $\cmat{N}{N}$. Similarly, for $\mathcal{X}$, a (closed) subspace of $\cmat{N}{N}$, we denote by $\meas_{\mathcal{X}}(\cdot)$ the Lebesgue measure on $\mathcal{X}$; thus, $\meas_{\mathcal{X}}(U)$ denotes the measure of $U$ as a subset of the space $\mathcal{X}$.

\noindent For $P\in\cmat{N}{N}$, by $\spr(P)$, the sparsity class induced by $P$, we will denote the subspace of matrices in $\cmat{N}{N}$ whose non-zero entry locations are a subset of those of $P$, i.e.,
\begin{align*}
\spr(P)=\left\{\acute{P}\in\cmat{N}{N}~:~P_{(n,l)}=0~\mbox{implies}~\acute{P}_{(n,l)}=0\right\}.
\end{align*}
Note that, by the above constructions, for $P,Q\in\cmat{N}{N}$, $\spr(P)$ is a subspace of $\spr(Q)$ if $\overline{P}\preceq\overline{Q}$.

\noindent\textbf{Problem formulation}. Suppose we have $N$ agents with (possibly complex-valued) data $z_{n}\in\mathbb{C}$, $i=1,\cdots,N$, distributed across the agents. Given a matrix $T\in\mathbb{C}^{N\times N}$, we are interested in computing the linear transformation $\mathbf{y}=T\mathbf{z}$, $\mathbf{z}=[z_{1}~z_{2}\cdots z_{N}]^{\top}$, in a distributed fashion using linear iterations of the form~\eqref{lin_dist_update} such that in a prescribed number of iterations, say $d$, each agent $n$ obtains the $n$-th component $y_{n}=\mathbf{t}_{n}^{\top}\mathbf{z}$ of the transform vector, with $\mathbf{t}_{n}^{\top}$ denoting the $n$-th row of $T$. To formalize, note that the distributed update at an iteration $k$ may be written in a compact form as 
\begin{align}
\label{lin_update_comp}
    \z(k)=A_{k}\z(k-1),
\end{align} 
where $\z(k)$ denotes the vector of the agent (computation states) and the weight matrix $A_{k}\in\spr(\oA)$, i.e., it respects the sparsity of the communication adjacency matrix $\oA$. The underlying design objective, i.e., to obtain a sequence of weight matrices to achieve finite-time computation, leads to the notion of $d$-step finite-time computability that we pursue in this paper:
\noindent\begin{definition}
\label{def:ftc}
A linear transformation $T\in\cmat{N}{N}$ is said to be finite-time computable in $d$-steps in a distributed fashion if there exist a sequence of matrices $A_{i}$, $i=1,\cdots,d$, in $\spr(\oA)$ such that, for all $\z\in\mathbb{C}^{N}$, the ($d$-step) iterative procedure
\begin{align*}
    \z(k)=A_{k}\z(k-1),~~k=1,\cdots,d,
\end{align*}
initialized with $\z(0)=\z$, yields $\z(d)=T\z$.
\end{definition}
Note that, in contrast to classical time-invariant consensus type protocols, we allow the weight matrices $A_{k}$'s to be time varying. In this paper, we aim to characterize the class of transforms $T$ that are $d$-step finite-time computable, and how to compute such transforms efficiently, i.e., in a small number of iterations $d$.

\noindent\textbf{A graph-structured matrix factorization based reformulation}. Clearly, we have the following equivalence between $d$-step finite-time computability of a transform $T$ and an associated graph-structured matrix factorization of $T$: the linear transform $T\in\cmat{N}{N}$ is finite-time computable in $d$-steps in a distributed fashion if and only if there exists a sequence of matrices $A_{i}$, $i=1,\cdots,d$, in $\spr(\oA)$ such that
\begin{align}
\label{graph_str_fact}
    T=A_{d}A_{d-1}\cdots A_{1}.
\end{align}
In what follows, we work with the factorization reformulation to analyze $d$-step finite computability.

\noindent\textbf{Key technical challenges}. We explain the key technical challenges through the matrix factorization based reformulation, i.e., obtaining $A_{k}$'s such that~\eqref{graph_str_fact} holds. Clearly, a major issue is that the factors $A_{i}\in\spr(\oA)$ could be sparse (if $\oA$ is sparse), whereas, the transform $T$ could be arbitrary, for instance, a dense matrix. We briefly review how existing approaches handle this issue in the rank-one case that is relatively well studied. For definiteness, consider a further special (but technically equivalent) case, the average consensus problem, i.e., $T=(1/N)\mathbf{1}\mathbf{1}^{T}$; in this case, the problem reduces to obtaining $A_{i}$'s, $i=1,\cdots,d$, such that $A_{i}\in\spr(\oA)$ and $A_{d}A_{d-1}\cdots A_{1}=(1/N)\mathbf{1}\mathbf{1}^{T}$. Usually, this is addressed in the setting of undirected connected graphs (i.e., symmetric and irreducible $\oA$) using eigenvalue matching techniques based on spectral properties of matrix polynomials. For instance, in~\cite{sandryhaila2014finite}, assuming $\oA$ is symmetric and irreducible, the idea is to design factors $A_{i}$'s of the form $A_{i}=\left(I+\alpha_{i}W\right)$, where $W\in\spr(\oA)$ is appropriately chosen (typically an irreducible symmetric doubly stochastic matrix). The resulting product is a polynomial in $W$ and the idea is to choose the $\alpha_{i}$'s such that all but one of the eigenvalues of the polynomial are nulled; further, it can be shown that this eigenvalue nulling requires at most $N$ steps or factors. However, this $O(N)$ matrix polynomial based idea breaks down in the case of general transforms $T$ as the underlying eigenvalue matching by using factors of the form $A_{i}=\left(I+\alpha_{i}W\right)$ becomes infeasible when $T$ is of higher rank. Nonetheless, a naive solution requiring $d=O(N^{2})$ iterations may be obtained for connected undirected graphs: we simply repeat $N$ rounds of computation, one for each of the linear combinations corresponding to each row of $T$ and noting that each of these linear combinations may be obtained in $O(N)$ iterations by leveraging  the above mentioned finite-time approaches.

Broadly speaking, in this work, we show that, under broad types of inter-agent connectivity conditions including strictly directed networks, the $O(N^{2})$ complexity can be reduced to $O(N)$ for almost all transformations. To this end, we develop a fundamentally different design approach based on ideas in algebraic geometry. Intuitively, the main idea is to cast the sparsity constrained matrix factorization formulation (see~\eqref{graph_str_fact}) as a system of polynomial equations in the \emph{free} components of the constituent factor matrices and obtain conditions on the structure (connectivity) of the communication adjacency matrix $\oA$ that guarantee solutions to~\eqref{graph_str_fact} exist with at most $d=O(N)$ factors.

\section{A Generic Structural Factorization}
\label{sec:fact} For $A\in\cmat{N}{N}$, denote by $\mathcal{E}(A)$ the set
\begin{align}
\label{def_calE}
\mathcal{E}(A)=\left\{(n,l)\in V^{2}~:~A_{nl}\neq 0\right\}.
\end{align}
Now, let $\overline{A},\overline{B}\in\{0,1\}^{N\times N}$, and suppose that $I_{N}\in\spr(\overline{A})$ and $I_{N}\in\spr(\overline{B})$, i.e., both $\overline{A}$ and $\overline{B}$ contain strictly non-zero diagonals. Consider the product $\overline{A}\overline{B}$ and the associated set $\mathcal{E}(\overline{A}\overline{B})$. Note,
\begin{align}
\label{def_calC}
\mathcal{E}(\overline{A}\overline{B})=\left\{(n,l)\in V^{2}~:~\sum_{m=1}^{N}\overline{A}_{(n,m)}\overline{B}_{(m,l)}\neq 0\right\}.
\end{align}
Since $\overline{A},\overline{B}$ have non-negative entries and both $\spr(\overline{A})$ and $\spr(\overline{B})$ contain the identity $I_{N}$, it readily follows that 
\begin{align}
\label{def_union}
\mathcal{E}(\overline{A})\cup\mathcal{E}(\overline{B})\subseteq\mathcal{E}(\overline{A}\overline{B}).
\end{align}
For each $(n,l)\in\mathcal{E}(\overline{A}\overline{B})$, let
\begin{align}
\label{def_union1}
\mathcal{E}_{(n,l)}(\overline{B})= & \left\{(m,l)\in\mathcal{E}(\overline{B})~:~\overline{A}_{(n,m)}\overline{B}_{(m,l)}\neq 0\right.\\ \nonumber & \left.~\mbox{for some $(n,m)\in\mathcal{E}(\overline{A})$}\right\},
\end{align}
and note, by construction, $\mathcal{E}_{(n,l)}(\overline{B})\neq\emptyset$ for all $(n,l)\in\mathcal{E}(\overline{A}\overline{B})$.

\noindent We make an observation here. The above development asserts that as long as $I_{N}\in\spr(\overline{A})$ and $I_{N}\in\spr(\overline{B})$, the sets of non-zero entries of $\overline{A}$ and $\overline{B}$ are both contained in the set of non-zero entries of the product $\overline{A}\overline{B}$ (see~\eqref{def_union}). In terms of associated digraphs, the above means that the edge-set of the digraph associated with the product matrix  $\overline{A}\overline{B}$ contains (at least) the edges in the digraphs associated with $\overline{A}$ and $\overline{B}$. In a sense, the product $\overline{A}\overline{B}$ leads to a digraph that is an \emph{extension} of the digraphs associated with $\overline{A}$ and $\overline{B}$. Intuitively, the above notion of extension will play a crucial role in obtaining a generic factorization of elements of $\cmat{N}{N}$ in terms of products of certain structured sparse matrices in $\cmat{N}{N}$. The following concept of simple extension will be useful.
\begin{definition}
\label{def:compatible} An ordered tuple $(\overline{A},\overline{B})$ of matrices in $\{0,1\}^{N\times N}$, with both $\overline{A}$ and $\overline{B}$ possessing strict non-zero diagonals, is said to be compatible if $\mathcal{E}(\overline{B})\subset\mathcal{E}(\overline{A})$.
\end{definition}
For a compatible tuple $(\overline{A},\overline{B})$, by $\beta_{(m,l)}(\overline{B})$, for each $(m,l)\in\mathcal{E}(\overline{B})$, we denote the set
\begin{align}
\label{def:compatible1}
\beta_{(m,l)}(\overline{B})=\left\{(n,l)\in\mathcal{E}(\overline{A}\overline{B})\setminus\mathcal{E}(\overline{A})~:~\overline{A}_{n,m}\overline{B}_{(m,l)}\neq 0\right\}.
\end{align}
In other words, for $(m,l)\in\mathcal{E}(\overline{B})$, $\beta_{(m,l)}(\overline{B})$ denotes the subset of new edges in $\mathcal{E}(\overline{A}\overline{B})$, i.e., those that are not present in $\mathcal{E}(\overline{A})$ (and hence in $\mathcal{E}(\overline{B})$), that make use of the edge $(m,l)$ in their generation. Note that $\beta_{(m,l)}(\overline{B})$ could be empty in general; clearly, $\beta_{(m,l)}(\overline{B})=\emptyset$ if $m=l$.
\begin{definition}
\label{def:simple_extension} Let $(\overline{A},\overline{B})$ be a compatible tuple of matrices in $\{0,1\}^{N\times N}$. We say that the product $\overline{A}\overline{B}$ is a simple extension of $(\overline{A},\overline{B})$ if $\left|\mathcal{\beta}_{(m,l)}(\overline{B})\right|\leq 1$ for all $(m,l)\in\mathcal{E}(\overline{B})$.
\end{definition}

\noindent Recall that, for $\overline{A}\in\BNN$, $\spr(\overline{A})$ is a (closed) subspace of $\cmat{N}{N}$. For $\overline{A},\overline{B}\in\BNN$ with strictly positive diagonal entries, we have $\acute{A}\acute{B}\in\spr(\overline{A}\overline{B})$ for all $\acute{A}\in\spr(\overline{A})$ and $\acute{B}\in\spr(\overline{B})$. Denote by $F_{\overline{A},\overline{B}}:\spr(\overline{A})\times\spr(\overline{B})\mapsto\spr(\overline{A}\overline{B})$ the product map $F_{\overline{A},\overline{B}}(\acute{A},\acute{B})=\acute{A}\acute{B}$ for $(\acute{A},\acute{B})\in\spr(\overline{A})\times \spr(\overline{B})$. Denote by $F_{\overline{A},\overline{B}}\left(\spr{\overline{A}}\times\spr(\overline{B})\right)\subseteq\spr(\overline{A}\overline{B})$ the range of this mapping.
 
\noindent We now obtain conditions under which the range $F_{\overline{A},\overline{B}}\left(\spr{\overline{A}}\times\spr(\overline{B})\right)$ is an everywhere dense subset of $\spr(\overline{A}\overline{B})$. Note, showing that $F_{\overline{A},\overline{B}}\left(\spr{\overline{A}}\times\spr(\overline{B})\right)$ is an everywhere dense subset of $\spr(\overline{A}\overline{B})$ is equivalent to showing that the system of alegbraic equations
\begin{align}
\label{lm:prod_dense2}
XY=\acute{C}
\end{align}
has a solution in $X\in\spr(\overline{A})$ and $Y\in\spr(\overline{B})$ for almost all $\acute{C}\in\spr(\overline{A}\overline{B})$. To this end, we have the following result.
\begin{theorem}
\label{lm:prod_dense} Let the ordered tuple $(\overline{A},\overline{B})$ of matrices in $\{0,1\}^{N\times N}$ be compatible in the sense of Definition~\ref{def:compatible} such that the product $\overline{A}\overline{B}$ is a simple extension of $(\overline{A},\overline{B})$ (in the sense of Definition~\ref{def:simple_extension}). Then, the range $F_{\overline{A},\overline{B}}\left(\spr{\overline{A}}\times\spr(\overline{B})\right)$ of the product mapping $F_{\overline{A},\overline{B}}(\cdot)$ contains an open dense subset of $\spr(\overline{A}\overline{B})$ such that
{\small
\begin{align}
\label{lm:prod_dense1}
\meas_{\spr(\overline{A}\overline{B})}\left(\cl\left(\spr(\overline{A}\overline{B})\setminus F_{\overline{A},\overline{B}}\left(\spr{\overline{A}}\times\spr(\overline{B})\right)\right)\right)=0,
\end{align}}
i.e., the system of alegbraic equations in~\eqref{lm:prod_dense2} has a solution in $X\in\spr(\overline{A})$ and $Y\in\spr(\overline{B})$ for almost all $\acute{C}\in\spr(\overline{A}\overline{B})$.
\end{theorem}
The proof of Theorem~\ref{lm:prod_dense} is omitted. First, we note that the algebraic system~\eqref{lm:prod_dense2} constitutes a system of polynomial equations in the \emph{free} entries of the matrices $X$ and $Y$, i.e., those entries $X_{n,l}$'s and $Y_{n,l}$'s that are unconstrained in $\spr(\overline{A})$ and $\spr(\overline{B})$ respectively. With this observation, the key technical step is to show that the system of equations is alegbraically independent~\cite{milne2012algebraic}, thus guaranteeing a solution for almost all $\acute{C}\in\spr(\overline{A}\overline{B})$; this involves a technical construction where we show that an associated Jacobian is of generic full rank by employing the structural assumption that the product $\overline{A}\overline{B}$ is a simple extension of $(\overline{A},\overline{B})$. 

\section{Main Results}
\label{sec:main _res} We now return to the problem of finite-time computability of a linear transform $T\in\cmat{N}{N}$ (in the sense of Definition~\ref{def:ftc}), or equivalently, the graph-structured matrix factorization reformulation~\eqref{graph_str_fact}. The results we state below on finite-time computability shows that, under appropriate connectivity conditions on the inter-agent communication adjacency matrix, almost all linear transforms $T$ may be computed in at most $N$ steps, i.e., in the sense of Definition~\ref{def:ftc}, almost all linear transforms are $N$-step finite-time computable.

\noindent Before proceeding to the case of general digraphs, we consider the case in which the adjacency matrix $\overline{A}$ corresponds to a directed cycle with self-loops on the $N$ vertices. We have the following result.
\begin{theorem}
\label{lm:fact_cycle}
Let the communication adjacency matrix $\overline{A}\in\BNN$ correspond to the directed cycle digraph with self-loops on the $N$ vertices, i.e., $\overline{A}_{n,n}=1$ for all $n\in V$, $\overline{A}_{n,n+1}=1$ for all $n\in V\setminus\{N\}$, $\overline{A}_{N,1}=1$, and all other entries of $\overline{A}$ are zero. Then, for almost all $T\in\cmat{N}{N}$, there exist matrices $A_{i}\in\spr(\overline{A})$, $i=1,\cdots,N$, such that
\begin{align}
\label{lm:fact_cycle1}
T=A_{N}A_{N-1}\cdots A_{1}.
\end{align}
\end{theorem}
The proof is omitted but essentially consists of a recursive application of Theorem~\ref{lm:prod_dense}. Roughly, by using the directed cycle structure of $\overline{A}$, we are able to verify the simple extension requirement of Theorem~\ref{lm:prod_dense} at each stage of the recursive construction, leading to the assertion that the set of matrices that may not be written as the product of $N$ factor matrices in $\spr(\overline{A})$ has a negligible (of Lebesgue measure zero) closure. To obtain this assertion and that we need at most $N$ factors, we also leverage the fact that the diameter of the directed cycle is bounded by $N$ and hence, intuitively, we are able to \emph{match} all the entries of $T$ (which could be a full matrix) with at most $N$ products. 

Theorem~\ref{lm:fact_cycle} may be extended to the case of general strongly connected digraphs as follows:
\begin{theorem}
\label{lm:fact_Hamilton}
Let the communication adjacency matrix $\overline{A}\in\BNN$ be such that $\overline{A}_{nn}=1$ for all $n=1,\cdots,N$ and there exists a Hamiltonian cycle in the associated digraph. Then, for almost all $T\in\cmat{N}{N}$, there exist matrices $A_{i}\in\spr(\overline{A})$, $i=1,\cdots,N$, such that $T=A_{N}A_{N-1}\cdots A_{1}$.
\end{theorem}
Note, the Hamiltonian cycle means that the associated communication digraph contains a directed cycle as a subgraph. The Hamiltonian cycle requirement is stronger than the strongly connectedness of the inter-agent communication digraph. By the Hamiltonian cycle assumption, by retaining the diagonal entries in $\overline{A}$ and by zeroing certain off-diagonal entries, we can obtain a matrix $\overline{A}_{\mbox{\scriptsize{red}}}\preceq\overline{A}$ in $\BNN$ that corresponds to a directed cycle with self-loops as used in the hypotheses of Theorem~\ref{lm:fact_cycle}. By restricting attention to factor matrices in $\spr(\overline{A}_{\mbox{\scriptsize{red}}})$ and noting that $\spr(\overline{A}_{\mbox{\scriptsize{red}}})\subseteq\spr(\overline{A})$, the desired factorization follows immediately as a consequence of Theorem~\ref{lm:fact_cycle}.

\section{Discussion}
\label{sec:disc} In this paper, we have considered finite-time (iterative) in-network computation of linear transforms via distributed linear iterations that conform to the sparsity of the (communication) network adjacency structure. By reformulating the finite-time computability problem to a sparsity constrained matrix factorization problem and employing tools from algebraic geometry, we have obtained a generic solution to the former problem: specifically, we have shown that almost all linear transforms are computable in at most $N$ steps (iterations) as long as the communication graph contains a directed Hamiltonian cycle. The results obtained extend the existing literature on finite-time distributed linear transform computation, that primarily focuses on the case of computing rank-one transforms (corresponding to homogeneous agent computation objectives) over undirected graphs, on multiple fronts including the computation of arbitrary transforms that correspond to potentially heterogeneous computation objectives at the network agents (i.e., each agent being interested in a possibly different linear combination of the network data) and over general directed communication graphs. 

There are several future research directions and extensions that are worth exploring. The above results, in a sense, provide a qualitative characterization, i.e., whether a transform $T$ is finite-time computable with guarantees on the computation time. However, in practice, the associated design question, i.e., how to obtain the factor matrices in~\eqref{lm:fact_cycle1} or equivalently the (time-varying) weights in the linear distributed iterations~\eqref{lin_dist_update}, needs to be handled in a systematic and efficient way. Given that $N$ factors are sufficient, we may resort to (non-convex) optimization techniques to obtain \emph{approximate} factorizations and hence weight designs. Secondly, the characterization may be conservative, i.e., for specific classes of transforms $T$ and depending on the communication digraph, fewer than $N$ factors (equivalently iterations) may be sufficient; this is a direction for future research. Another research issue, related to the above, is to explore to what extent the Hamiltonian cycle requirement may be relaxed or modified.

\bibliographystyle{IEEEtran}
\bibliography{finite_time_comp.bib}

\end{document}